\documentclass[11pt]{amsart}
\usepackage{amsmath, amssymb, amsfonts, amscd}
\newcommand{\proj}[1]{\operatorname{Proj}(#1)}

\newcommand{\length}{\operatorname{length}}

\newcommand{\rank}{\operatorname{rank}}

\newcommand{\A}{\operatorname{A}}

\newcommand{\ch}{\operatorname{ch}}
\newcommand{\K}{\operatorname{K_0}} 
\newcommand{\G}{\operatorname{G_0}} 

\newtheorem{theorem}{Theorem}
\newtheorem{prop}{Proposition}

\newtheorem{lemma}{Lemma}

\title{\sc  A Correspondence
Between Hilbert Polynomials and Chern Polynomials over Projective
Spaces }

\author{\sc C-Y. Jean Chan }

 
\begin{document}

\begin{abstract}
  
  We construct a map $\zeta$ from $\K(\mathbb P^d)$ to $(\mathbb
  Z[x]/x^{d+1})^{\times} \times \mathbb Z$, where $(\mathbb
  Z[x]/x^{d+1})^{\times}$ is a multiplicative Abelian group with
  identity $1$, and show that $\zeta$ induces an isomorphism between
  $\K(\mathbb P^d)$ and its image. This is inspired by a
  correspondence between Chern and Hilbert polynomials stated in
  Eisenbud~\cite[Exercise~19.18]{E}. The equivalence relation of these
  two polynomials over $\mathbb P^d$ is discussed in this paper. 
  
\end{abstract}

\maketitle

\section{Introduction}~\label{intro}

It is known that the Hilbert polynomial of a coherent sheaf over a
projective scheme is closely related to the Chern polynomial of the
sheaf by Hirzebruch-Riemann-Roch theorem. Throughout the paper,
$\mathbb P^d$ denotes a projective space over an algebraically closed
field. In fact, over $\mathbb P^{d}$, knowing the
Hilbert polynomial is equivalent to knowing the Chern polynomial. This
fact is pointed out in Eisenbud~\cite[Exercise~19.18]{E}. We will
briefly describe these definitions associated with a coherent sheaf in
the next section.  The Chern and Hilbert polynomials are
quite different in terms of degrees and coefficients. The
Hirzebruch-Riemann-Roch theorem makes a connection from one to the
other. In this paper, we let $\mathcal A_0$ and $\mathcal B$ denote
two Abelian groups which are generated by Chern polynomials and
Hilbert polynomials respectively.  We prove in Theorem~\ref{thmisom}
the existence of an isomorphism between the Grothendieck group
$\K(\mathbb P^d)$ and $\mathcal A_0 \times \mathbb Z$.  This appears
to be an analogous fact that $\K(\mathbb P^d)$ and $\mathcal B$ are
isomorphic as it is shown in \cite{E}.  Let $P(t)$ and $C(x)$ denote
the Hilbert and Chern polynomials of a coherent sheaf respectively.
Our main discussions on the isomorphisms mentioned in the above
conclude the following three equivalent statements: for any two
coherent sheaves $\mathcal M$ and $\mathcal N$, \\ \indent (1)
$\mathcal M$ and $\mathcal N$ represent the same class in $\K(\mathbb
P^d)$, \\ \indent (2) $P_{\mathcal M}(t) = P_{\mathcal N}(t)$, \\
\indent (3) $C_{\mathcal M}(x)=C_{\mathcal N}(x)$ and $\rank \mathcal
M= \rank \mathcal N$. \\ It is easy to see (1) $\Rightarrow$ (2) and
(1) $\Rightarrow$ (3) by definitions.  (2) $\Rightarrow$ (1) is done
in \cite{E}. We give a proof for (3) $\Rightarrow$ (1).

The paper is arranged in the following way. Section~\ref{background}
contains some necessary background materials.  It is not possible to
give precise definitions in this paper. We state their important
properties which are often used in our discussion. Section~\ref{isoms}
describes two isomorphic group structures on $\K(\mathbb P^d)$. One is
induced by a map $\eta: \K(\mathbb P^d) \longrightarrow \mathcal B$,
$\mathcal B \subset (\mathbb Q[t]/(t^{d+1}))^{+} $
(\cite[Exercises~19.16 and~19.17]{E}) and the other is induced by 
$\zeta: \K(\mathbb P^d) \longrightarrow \mathcal A$, $\mathcal A
\subset (\mathbb Z[x]/(x^{d+1}))^{\times} \times \mathbb
Z$~(Theorem~\ref{thmisom}).  The above equivalent statements on
equivalent classes, $P_{\mathcal M}(t)$ and $C_{\mathcal M}(x)$ then
follow from Theorem~\ref{thmisom}. The isomorphism from $\mathcal A$
to $\mathcal B$ induced by $\zeta^{-1}$ and $\eta$ recovers the
Hirzebruch-Riemann-Roch theorem. This result together with the fact
provided in Eisenbud~\cite[Exercise~19.18]{E} gives an algorithm which
computes the Hilbert polynomial of a coherent sheaf from its Chern
polynomial and vise versa without using an explicit sheaf structure.
This one-to-one correspondence between the two polynomials is
discussed in Section~\ref{equivsection}.

\section{Background}\label{background}

We begin with the definitions of some notations which are used
throughout the paper.  Let $\mathbb P^d$ be a projective space over an
algebraically closed field. For any coherent sheaf $\mathcal M$ over
$\mathbb P^d$, $\mathcal M$ is associated with a graded finitely
generated module $M=\oplus_n M_n$ over the polynomial ring
$S=k[x_0,\dots, x_d]$ with homogeneous components $M_n$. Conversely, a
graded module also defines a coherent sheaf, but there are usually
more than one graded module associated to a given coherent sheaf ({\em
  c.f.}  \cite{Hart, K}). For each $\mathcal M$, there exists the {\em
    Hilbert polynomial} $P_{\mathcal M}(t)$ of $\mathcal M$ such that
  for any large enough integer $n \in \mathbb N$, $P_{\mathcal M}(n)$
  coincides the value of the Hilbert function of the module $M$,
  $H_M(n)=\length(M_n)$.  For example, any twisted structure sheaf
  $\mathcal O(-m)$ with $m \in \mathbb Z$ is associated with the
  graded module $S[-m]$ and $P_{\mathcal O(-m)}(t)= \binom{t+d-m}{d}$.
  The Hilbert polynomials are additive on short exact sequences of
  sheaves.

The {\em Chern polynomial} $C_{\mathcal M}(x)$ (often called the {\em
  total Chern class}) of $\mathcal M$ is a formal sum of the Chern classes
which, in geometry, are usually viewed as cycles in the cohomology groups
({\em c.f.} Griffiths and Harris~\cite{GH}) or operators on the Chow
groups~({\em c.f.} Fulton~\cite{F} or Roberts~\cite{R1}). The Chern
classes considered in this paper are in the latter form.  In general,
their definitions are very complicated.  We describe these notions for
sheaves over $\mathbb P^d$ and recall a few properties that will be
useful for the latter discussions. The complete details can be found
in the above cited references.

The {\em Chow group} $\A_*(\mathbb P^d)$ of $\mathbb P^d$ is generated
by the linear subspaces $\mathbb P^{d-\ell}, \ell=0,\dots, d$, so it
has a simple structure; $\A_*(\mathbb P^d) \cong \mathbb
Z^{d+1}$. Therefore, the Chern classes can be identified with
integers. For a locally free sheaf $\mathcal M$ of finite rank $r$, there
exist $r$ Chern classes $c_1(\mathcal M), \dots, c_r(M)$, and the
{\em Chern polynomial} of $\mathcal M$ is defined to be a formal sum
of $c_i(\mathcal M)$,
\[C_{\mathcal M}(x)= 1+ \sum_{i=1}^r c_i(\mathcal M)x^i=
1+c_1(\mathcal M)x + \cdots + c_r(\mathcal M)x^r \pmod{x^{d+1}}.\] For
instance, $c_1(\mathcal O(-m))=-m$ for any twisted sheaf with $m \in
\mathbb Z$ and there is no higher Chern classes. Thus, $C_{\mathcal
  O(-m)}(x)=1-mx$.  An important property called the {\em Whitney sum
  formula} states that the Chern polynomials are multiplicative on 
short exact sequences of sheaves.

Over a nonsingular variety, every coherent sheaf admits a unique
minimal resolution of locally free sheaves up to quasi-isomorphisms
\begin{equation}\label{resolution}
 0 \rightarrow \oplus_{j_d} \mathcal O(-j_d)^{\beta_{d,j_d}} \rightarrow \dots
 \rightarrow \oplus_{j_1} \mathcal O(-j_1)^{\beta_{1,j_1}} \rightarrow
   \oplus_{j_0} \mathcal O(-j_0)^{\beta_{0,j_0}} \rightarrow \mathcal M
\rightarrow 0.
\end{equation}
Using the Whitney sum formula, the definition of Chern classes
can be extended to the coherent sheaves.
By (\ref{resolution}), the Chern polynomial of
$\mathcal M$ over $\mathbb P^d$ is a polynomial modulo $x^{d+1}$ with 
integer coefficients
\begin{equation}\label{eqchern}
  C_M(x) = \frac {{\displaystyle \prod _{i:\mbox{even}} \prod_{j_i} 
  (1-j_i x)^{\beta_{i,j_i}} }}
  {{\displaystyle \prod_{i:\mbox{odd}} \prod_{j_i} 
   (1-j_i x)^{\beta_{i,j_i}} }}
  \mbox{\hspace{.3in}({\rm mod} $x^{d+1}$)} 
\end{equation}
and the Hilbert polynomial is of degree at most $d$ with rational
number coefficients
\[P_M(t) \begin{array}[t]{cl}
  =& { \sum_{i,j_i}(-1)^{\beta_{i,j_i}}P_{S(-j_i)}(t)}.
 \end{array} \]
 
 Recall that $C_{\mathcal O(-m)}(x)=1-mx$ and $P_{\mathcal
   O(-m)}(t)=\binom{t+d-m}{d}$. The Chern polynomials of locally free
 sheaves always have integer coefficients and the degree varies while
 the Hilbert polynomials of such sheaves have rational number
 coefficients and the degree is fixed by $\dim X=d$.  A consequence of
 the Hirzebruch-Riemann-Roch theorem shows the connection of the Euler
 characteristic and Chern characters of coherent sheaves which leads a
 representation of Hilbert polynomial in terms of Chern classes in
 some special cases. We will recall this in Section~\ref{equivsection}.
 
 Next, we define the Grothendieck group. The Grothendieck group of
 locally free sheaves, denoted $\K(\mathbb P^d)$, is the Abelian group
 generated by all the locally free sheaves $[\mathcal M]$ modulo the
 subgroup generated by $[\mathcal M] - [\mathcal M'] - [\mathcal M'']$
 whenever $0 \rightarrow \mathcal M' \rightarrow \mathcal M
 \rightarrow \mathcal M'' \rightarrow 0$ forms an exact sequence of
 sheaves. We also denote $[\mathcal M]$ by $[M]$ if the sheaf
 $\mathcal M$ is associated with the module $M$. The Grothendieck
 group of coherent sheaves $\G(\mathbb P^d)$ is defined in a similar
 way. Since every coherent sheaf admits a locally free resolution in
 the form of (\ref{resolution}), $\K(\mathbb P^d)$ is isomorphic to
 $\G(\mathbb P^d)$. Henceforth, we use $\K(\mathbb P^d)$ and refer it
 as the Grothendieck group of $\mathbb P^d$ for simplicity.  The
 generators of $\K(\mathbb P^d)$ can be described precisely in the
 followings.  The module $S/(x_0,\dots, x_d)$ defines a zero sheaf
 since $(x_0,\dots,x_d)$ is an irrelevant ideal. We take the Koszul
 resolution of $S/(x_0,\dots, x_d)$ and obtain a long exact sequence
 on locally free sheaves,
\begin{equation}\label{sqsheaves}
0\rightarrow \mathcal O(-d-1) \rightarrow
(\mathcal O(-d))^{\binom{d+1}{d}}
 \rightarrow \cdots \rightarrow
 (\mathcal O(-1))^{\binom{d+1}{1}} \rightarrow \mathcal O
 \rightarrow 0. \end{equation} 
This gives a relation for the twisted
 sheaves in $\K(\mathbb P^d)$ which expresses $[\mathcal O(-d-1)]$ as
 an alternating sum of $[\mathcal O], [\mathcal O(-1)], \dots, 
 [\mathcal O(-d)]$. If we twist the exact sequence~(\ref{sqsheaves}) by 
 a degree, say by
 degree one, then we have the following exact sequence \[ 0\rightarrow
 \mathcal O(-d) \rightarrow \mathcal (\mathcal
 O(-d+1))^{\binom{d+1}{d}} \rightarrow \cdots \rightarrow (\mathcal
 O)^{\binom{d+1}{1}} \rightarrow \mathcal O(1) \rightarrow 0. \] Thus,
 $[\mathcal O(1)]$ is also generated by the same set of twisted
 sheaves and similarly for other degrees. This implies that
 $\K(\mathbb P^d)$ is generated by $[\mathcal O (-m)]$ with
 $m=0,\dots,d$. Furthermore, these are free generators. A
 brief argument for this fact will be developed in the next
 section. On the other hand, let $S_{\ell}$ denote the graded module of 
$S$ modulo $\ell$ variables
\[ S_{\ell} = k[x_0, \cdots, x_d]/ ( x_{d-\ell +1}, \cdots, x_d). \] 
The Koszul complex on $x_{d-\ell+1},\dots, x_d$ provides a
 resolution of locally free sheaves for $S_{\ell}$. By a standard
 argument, $[S_{\ell}]$, $\ell=0,\dots,d$, also generates
 $\K(\mathbb P^d)$.

 It is known that different sheaves may have the same Hilbert
 polynomials and Chern polynomials. However, both polynomials are
 well-defined for the equivalence classes of coherent sheaves in the
 Grothendieck group. We have already
 seen in Eisenbud~\cite{E} that the Hilbert polynomials characterize
 the classes in $\K(\mathbb P^d)$.  The main goal of this paper is to
 show that the Chern polynomials do the same job; namely, distinct classes
have different pairs of Chern polynomial and rank.

\section{Groups Isomorphic to $\K(\mathbb P^d)$}\label{isoms}

If a polynomial with rational coefficients has integral values at
large integers, then it can be written as a linear combination over
$\mathbb Z$ of the following binomial coefficient functions in $t$
\[ \binom{t}{0}, \binom{t}{1}, \binom{t}{2}, \dots, \binom{t}{\ell}, \dots .\]
These polynomials can be replaced by
\[ \binom{t}{0}, \binom{t+1}{1}, \binom{t+2}{2}, \dots,
\binom{t+\ell}{\ell}, \dots . \]
If we let $a_{\ell }=\binom{t}{\ell}$ and
$b_{\ell}=\binom{t+\ell}{\ell}$, then $b_{\ell} =
\sum_{i=0}^{\ell} \binom{\ell}{i}a_i$ and $a_{\ell} =
\sum_{i=0}^{\ell}(-1)^i \binom{\ell}{i}b_i$.

Let $\mathbb P^d=\proj S=\proj{k[x_0, \cdots, x_d]}$ be as in the
previous section.  Then, $ \binom{t+d}{d}$ is exactly the Hilbert
polynomial of $\mathcal O_{\mathbb P^d}$ and $P_{s_{\ell}}(t)=
\binom{t+d-\ell}{d-\ell}$. Since Hilbert polynomials have integral
values at large integers, for any graded module $M$,
$P_{M}(t)$ can be written as a linear combination of $P_{S_{\ell}}(t),
\ell \in \mathbb N \cup \{0\}$.  Let $\mathcal B$ denote the Abelian
group generated by all Hilbert polynomials of coherent sheaves over
$\mathbb P^d$.  The group $\mathcal B$ is a subgroup of the additive
group $(\mathbb Q[t]/t^{d+1})^+$ with identity $0$. Then,
$P_{S_{\ell}}(t)$, $\ell =0,\dots,d$, form a set of generators for
$\mathcal B$. Moreover, these generators are linearly independent
since $\deg P_{S_{\ell}}(t)= \ell$. Thus,
$\mathcal B$ is a free Abelian group of rank $d+1$.

Let $\alpha$ be a class in $\K(\mathbb P^d)$ represented by some sheaf
$\mathcal M$. The map \[\eta : \K(\mathbb P^d) \longrightarrow
\mathcal B \] takes $\alpha$ to the Hilbert polynomial $P_{\mathcal
  M}(t)$ of $\mathcal M$ induces an isomorphism. To see this, we note
that $\eta$ is surjective since $P_{S_{\ell}}$, $\ell =0, \dots, d$,
generate $\mathcal B$. That these generators are linearly independent
implies $[S_{\ell}]$ are also linearly independent. Therefore,
$\K(\mathbb P^d)$ is generated freely by $[S_{\ell}]$, $\ell=0, \dots,
d$, and the injectivity follows ({\em cf.}~\cite[Exercise~19.17]{E}).
Another proof, using $\{[\mathcal O(-m)]: m=0, \dots, d\}$ as a
generating set and in which $\eta$ is induced by a map taking
$[\mathcal O(-m)]$ to its Hilbert series, can be found also
in~\cite[Exercise~19.16]{E} in great details.

Alternatively, let $\mathcal A_0$ denote the Abelian group generated
by all the Chern polynomials of the coherent sheaves over $\mathbb
P^d$. Similar to $\mathcal B$, $\mathcal A_0$ is a subgroup of the
Abelian multiplicative group $(\mathbb Z[x]/x^{d+1})^{\times}$ with
identity $1$. Let $\mathcal A$ denote the subgroup $\mathcal A_0
\times \mathbb Z$ of $(\mathbb Z[t]/t^{d+1})^{\times} \times \mathbb
Z$ which has the natural group structure that for any two elements
$(f(x),r)$ and $(g(x),s)$, $(f(x),r)+(g(x),s)=(f(x)g(x),r+s)$ and
$(1,0)$ is the identity. For any $\alpha$ in $\K(\mathbb P^d)$
represented by a locally free sheaf $\mathcal M$, we define a map from
$\K(\mathbb P^d)$ to $\mathcal A = \mathcal A_0 \times \mathbb Z$, \[
\begin{array}{cccl} \zeta:& \K(\mathbb P^d) &\longrightarrow &
\mathcal A =
 \mathcal A_0 \times \mathbb Z \\ 
 & \alpha=[\mathcal M] & \longrightarrow & 
 (C_{\mathcal M}(x), \rank \mathcal M). \end{array} \] 
 The map $\zeta$ is a well-defined group homomorphism by the
 Whitney sum formula. It should be noted that the component
 $\mathbb Z$ in $\mathcal A$ is necessary in order to distinguish
 different classes which have the same Chern polynomial.  The
 simplest examples are $\alpha=[\mathcal O_{\mathbb P^d}]$ of rank one
 and $\displaystyle{\beta=[ \oplus_r\mathcal O_{\mathbb P^d}]}$ of
 rank $r \neq 0,1$.  Both $\alpha$ and $\beta$ have the Chern
 polynomial equal to 1 while $\beta = r \alpha \neq \alpha$ in
 $\K(\mathbb P^d)$. Analogous to the isomorphism defined by $\eta$, we
 prove that $\zeta$ is also an isomorphism.

\begin{theorem}\label{thmisom}
$\zeta: \K(\mathbb P^d) \longrightarrow \mathcal A$ is an isomorphism
of Abelian groups.  
\end{theorem}

The following Lemma~\ref{freegen} implies that $\mathcal A$ is free of
rank $d+1$. The generators $(1-\ell x, 1)$ of $\mathcal A$ are the
image of $[\mathcal O(-\ell)]$ for all $\ell$ so $\eta$ is surjective
and therefore, it is an isomorphism since both groups are free of the
same rank.

\begin{lemma}\label{freegen}
The group $\mathcal A_0$ is freely generated by $1-x, \dots,
1-dx$. Furthermore, $(1,1), (1-x,1), \dots,(1-dx,1)$ are free
generators for $\mathcal A$.
\end{lemma}

\begin{proof}
  It is clear that $\mathcal A_0$ is generated by $1-x,\dots,1-dx$ and
  that $\mathcal A$ is generated by $(1,1), (1-x,1), \dots, (1-dx,1)$
  by the resolutions~(\ref{resolution}) and (\ref{sqsheaves}).
If   
\[ r_0(1,1) +
r_1(1-x,1) + \cdots + r_d(1-dx, 1)=(1,0) \] in $\mathcal A$ for some
$r_0, \dots, r_d \in \mathbb Z$.  Then, the following (\ref{eq1})
and~(\ref{eq2}) hold,
\begin{eqnarray}
&&(1-x)^{r_1}\cdots(1-dx)^{r_d} \equiv 1
 \pmod{ x^{d+1}}, \label{eq1} \\
&& r_0+r_1 + \cdots + r_d =0.  \label{eq2} \end{eqnarray}
It suffices to show that $ 1-x,
  \dots, 1-dx$ are linearly independent; that is, (\ref{eq1}) implies
  $r_1=\cdots=r_d=0$. Then, the linearly independence of $(1,1),
  (1-x,1), \dots, (1-dx,1)$ follows from (\ref{eq2}).

Without loss of generality, we do the following argument assuming that
none of $r_1, \dots, r_d$ is zero. (If any of $r_0, \dots, r_d$ is
zero, then it follows a similar argument which lead to the same
contradiction.) We take the derivative of the equation in (\ref{eq1})
and obtain 
\[(1-x)^{r_1}(1-2x)^{r_2}\cdots(1-dx)^{r_d}
 \left( \frac{-r_1}{1-x}+\frac{-r_2}{1-2x} + \cdots + \frac{-r_d}{1-dx} \right)
 \equiv 0 (\operatorname{mod} x^{d}).\]
The above product is taken in the unique factorization domain
$\mathbb Z[[x]]$ and a simple computation shows that
\[ (1-x)^{r_1}(1-2x)^{r_2}\cdots(1-dx)^{r_d}=1-(r_1+2r_2+\cdots +
  dr_d)x + \cdots \not\equiv 0 (\operatorname{mod}x^{d}). \]
Therefore, 
\[\frac{r_1}{1-x}+\frac{2r_2}{1-2x} + \cdots + \frac{dr_d}{1-dx} \equiv 0
 (\operatorname{mod} x^{d}). \]
Using the Taylor expansion, we have
\begin{equation}\label{equiv}
\begin{array}[t]{cc}
(r_1+2r_2+\cdots+dr_d)+(r_1+2^2r_2+\cdots+d^2r_d)x+\cdots \\ \\
 + (r_1+2^dr_2+ \cdots + d^dr_d)x^{d-1}
    \equiv 0 (\operatorname{mod} x^{d}). \end{array} \end{equation}
The equivalence given by Equation~(\ref{equiv}) provides a linear
system in $r_1, \dots, r_d$ with Vandemonde coefficients if
$r_1,\dots,r_d$ are all nonzero. This is a contradiction because a
Vandemonde system has only trivial solutions.  Therefore,
$r_1=r_2=\cdots=r_d=0$ and $r_0=0$ by~(\ref{eq2}).  This complete the
proof of both assertions in the lemma.
\end{proof}

The isomorphisms $\eta$ and $\zeta$ shows that the three groups
$\K(\mathbb P^d)$, $\mathcal A$ and $\mathcal B$ are 
isomorphic and the following conditions.

\begin{theorem}\label{tfae} 
For any coherent sheaves $\mathcal M$ and $\mathcal N$ on $\mathbb
P^d$, the followings are equivalent: \\
\indent $(1)$ $[\mathcal M]$ = $[\mathcal N]$ in $\K(\mathbb P^d)$. \\
\indent $(2)$ $P_{\mathcal M}(t) =  P_{\mathcal N}(t)$. \\
\indent $(3)$ $C_{\mathcal M}(x) = C_{\mathcal N}(x)$ and $\rank
\mathcal M = \rank \mathcal N$. 
\end{theorem}

\section{The equivalence of Chern and Hilbert polynomials}\label{equivsection}

This section discusses the close relationship of the Chern and
Hilbert polynomials which inspires the work presented in the previous
section. The Hirzebruch-Riemann-Roch theorem
relates the Euler characteristic with the Chern characters.
Not much is known about representing the Hilbert
polynomial in terms of Chern classes in general. We will discuss this and its
converse over $\mathbb P^d$ in the current section. 
Proposition~\ref{chernhilbert} provides a one-to-one correspondence
between the two polynomials and an algorithm for a computational purpose. 
 
The Hirzebruch-Riemann-Roch theorem proves that there exists a certain
maps from the Grothendieck group of a scheme $X$ to its Chow group and
that it commutes with the maps induced by a projective map from $X$ to
a point. This theorem for projective spaces induces an expression of
the Hilbert function in terms of the Chern classes. In order to
make it precise, we need to introduce the Chern characters of $\mathcal
M$. Suppose the Chern polynomial can be decomposed into
\begin{equation}\label{eqroots}
C_{\mathcal M}(x)=(1-\alpha_1 x) \cdots (1- \alpha_d x) \end{equation}
in $(z[x]/(x^{d+1})^{d+1})$. 
In this case, $\alpha_1, \dots, \alpha_d$ are called {\em Chern
  roots}. 
Then, the {\em Chern character} of $\mathcal M$ is a power series defined by
\begin{equation}\label{eqchar}
\ch(\mathcal M)=e^{\alpha_1 x}+ \cdots + e^{\alpha_d x}. \end{equation}
The coefficient of the $x^i$ in the Taylor expansion of (\ref{eqchar})
is called the {\em i}-{\em th Chern character} of $\mathcal M$, denoted 
$\ch_i(\mathcal M)$. Since each
$\ch_i(\mathcal M)$ is a symmetric function of $\alpha_i$ and the Chern
classes are elementary symmetric functions in $\alpha_i$, the Chern
characters $\ch_i(\mathcal M)$ can be expressed as a polynomial in the
Chern classes. The first few terms are
\begin{equation}\label{chernchar} \begin{array}[t]{ccl}
\ch(\mathcal M)& = & r + c_1 x + \frac{1}{2!}(c_1^2-2c_2)x^2 +
  \frac{1}{3!}(c_1^3-3c_1c_2+3c_3)x^3 + \\ \\
 &&  \frac{1}{4!}(c_1^4 -4c_1^2c_2+4c_1c_3+2c_2^2-4c_4)x^4 + \cdots,
\end{array} \end{equation}
where $c_i=c_i(\mathcal M)$ and $r=\rank \mathcal M$
(see~\cite[Example~15.1.2]{F} for the exact formulations).
We should note that the factorization~(\ref{eqroots}) does not always
exit over the current projective scheme. However,
the expressions in~(\ref{chernchar}) are independent from the existence 
of the Chern roots.

For any power series $s(x)$ in $x$, let $\Phi(s(x))$ denote the {\em
  coefficient of $x^d$ in the Taylor expansion of the expression
  $\displaystyle{s(x)\left(\frac{x}{1-e^{-x}}\right)^{d+1}}$}.
Theorem~\ref{thmHRR} states the Hirzebruch Riemann-Roch theorem for
$\mathbb P^d$.  The details can be found in the references by Fulton
and Lang~(\cite[Example~15.1.4]{F},~\cite{FL}) and
Hirzebruch~(\cite[Lemma~1.7.1]{Hirz}).
\begin{theorem}[Hirzebruch-Riemann-Roch]\label{thmHRR}
Let $X=\mathbb P^d$. Then, for any locally free sheaf $\mathcal M$ on $X$,
\begin{equation}\label{eqHRR}
\Phi(\ch(\mathcal M))
= \chi(\mathcal M),
\end{equation}
where $\chi(\mathcal M)=\sum_{i \geq 0}(-1)^i \dim_k H^i(X, \mathcal M)$, 
the alternating sum of the cohomology groups, is the Euler characteristic
of $\mathcal M$.
\end{theorem}

Let $\mathcal M$ be defined by a finitely generated graded module
$M=\oplus M_n$.  By the induction on dimension of the support of
$\mathcal M$, it shows that $H^i(X, \mathcal M(n))=0$ for all $i>0$
and
\[\chi(\mathcal M(n))=\sum_{i \geq 0} (-1)^i \dim H^i(X, \mathcal M(n))=
\dim H^0(X,\mathcal M(n))= \dim M_n, \] for $n \gg 0$ ({\em cf.}
Hartshorne~\cite[Chapter~III]{Hart}).  In the case where $\mathcal M$
is locally free, we have
\begin{equation}\label{HRRproj}
 \Phi(\ch(\mathcal M(n))) = \dim H^0(X,\mathcal M(n)) = \dim M_n
\end{equation}
by the Hirzebruch-Riemann-Roch theorem.  In particular, $ \ch(\mathcal
M(n))=\ch(\mathcal M \otimes \mathcal O(n)) = \ch(\mathcal
M)\ch(\mathcal O(n))= e^{nx}\ch(\mathcal M)$. We replace $n$ in
(\ref{HRRproj}) by an intermediate $t$.  The left hand side of
(\ref{HRRproj}), $\Phi(e^{tx}\ch(\mathcal M))$, becomes a polynomial
in $t$ whose values at large integers agree with values of the Hilbert
function of $M$. Therefore, it is the Hilbert polynomial ({\em cf.}
Fulton~\cite[Example~15.2.7(a)]{F}).

If $C_{\mathcal M}(x)=1+c_1(\mathcal M)x+ \cdots +
c_r(\mathcal M)x^r$ is the Chern polynomial of some locally free sheaf
$\mathcal M$ of rank $r$, then $\ch(\mathcal M)$ is known explicitly
as it is shown in (\ref{chernchar}).
The Hilbert polynomial $P_{\mathcal M}(t)$ of
$\mathcal M$ obtained by $\Phi(e^{tx}\ch(\mathcal M))$ is
the coefficient of $x^d$ in $e^{tx}\ch(\mathcal M)
 \left(\frac{x}{1-e^{-x}}\right)^{d+1}$. Part~A of the following
 proposition is an immediate consequence of Theorem~\ref{thmHRR} as it
 is explained in the above.

\begin{prop}~\label{chernhilbert}
Let $X=\mathbb P^d$ and let $\mathcal M$ be a coherent sheaf on
$X$ of rank $r$.

{\bf A.} Let $C_{\mathcal M}(x)$ be the Chern polynomial of $\mathcal M$.
      Then the Hilbert polynomial of $\mathcal M$ is
\begin{equation}\label{eqCP} P_{\mathcal M}(t)=
  \Phi(e^{tx}\ch(\mathcal M))  .\end{equation}

{\bf B.} If $\displaystyle{P_{\mathcal M}(t)= \sum_{\ell=0}^d
      a_{\ell} \binom{t+d-\ell}{d-\ell}}$, then
      \begin{equation}\label{eqPC} C_{\mathcal M}(x) \equiv \prod_{\ell=0}^d
         [C_{S_{\ell}}(x)]^{a_{\ell}}, \pmod{x^{d+1}},
          \end{equation}
where $C_{S_{\ell}}(x)$ is the Chern polynomial of $S_{\ell}$ as a module 
over $S$.
\end{prop}

\begin{proof} It remains to prove Part~B. We recall from
  Section~\ref{background} that $S_{\ell}$ defines the linear subspace
  $\mathbb P^{d-\ell}$ in $\mathbb P^d$. Since the Hilbert polynomial
  of $S_{\ell}$ is $\binom{t+d-\ell}{d-\ell}$, by the assumption and
  Theorem~\ref{tfae}, $\displaystyle{P_{\mathcal M}(t)=
    \sum_{\ell=0}^d a_{\ell} \binom{t+d-\ell}{d-\ell}=\sum_{\ell=0}^d
    a_{\ell} P_{S_{\ell}}(t)} $ in $\mathcal B$ if and only if
  $[\mathcal M]=\sum_{\ell=0}^d a_{\ell} [S_{\ell}]$ in
  $\K(\mathbb P^d)$. Therefore, $\displaystyle{C_{\mathcal M}(x) \equiv
    \prod_{\ell=0}^d (C_{S_{\ell}}(x))^{a_{\ell}}} \pmod{x^{d+1}}$.
\end{proof} 

We note that the expression of $P_{\mathcal M}(t)$ in the hypothesis
of Part~B is always possible since $\binom{t+d-\ell}{d-\ell}, \ell=0,
\dots, d$, form a
basis for the group $\mathcal B$ of all the Hilbert polynomials.  The
correspondence between the two polynomials can be viewed explicitly in
the followings. Let $\sigma$ denote the group homomorphism $\eta \circ
{\zeta}^{-1}$
\[ \sigma=\eta \circ {\zeta}^{-1} : 
\mathcal A \longrightarrow \mathcal B. \] 
Any element in the group $\mathcal A$ of the Chern polynomials
can be written as $\displaystyle{(\prod_{m=1}^d(1-mx)^{r_m},s)}$ for
some $r_m$ and $s$ in $\mathbb Z$.  It is not hard to see that such an
element has a preimage via $\zeta^{-1}$ in $\K(\mathbb P^d)$ as
$\sum_{m=1}^d a_m [\mathcal O(-m)] + (s-r)[\mathcal O]$ where $a=a_1
+\cdots +a_m$.  This implies
\[\begin{array}{ll}
 \sigma((\prod(1-mx)^{a_m},s))& =\sum_{m=1}^{d} a_m \binom{t+d-m}{d}
+(s-a)\binom{t+d}{d} \\ \\ 
 & = \sum_{m=1}^d a_m P_{\mathcal O(-m)}(t) + (s-a) P_{\mathcal
   O}(t). \end{array}\] 
An element $(f(x), s)$ in $\mathcal A$ is said
to be {\em representative by $\mathcal M$ } (or {\em $\mathcal M$
  represents $(f(x), s)$}) if there exists a sheaf $\mathcal M$ such
that $f(x)$ is the Chern polynomial of $\mathcal M$ and $\rank
\mathcal M=s$.  Thus, for any representative $(f(x),s)$ in $\mathcal A$, 
$\sigma$ takes $(f(x),s)$ to the
Hilbert polynomial of $M$. The computation can be carried out by
(\ref{eqCP}). Conversely, the preimage of the Hilbert
polynomial of $\mathcal M$ is the pair of the Chern polynomial and the
rank of $\mathcal M$. This preimage is uniquely determined since
$\sigma=\eta \circ \zeta^{-1}$ is an isomorphism.  
More precisely, (\ref{eqPC}) computes the Chern polynomial and $a_0 + \cdots
+ a_d$ indicates the rank which is independent from the choices of a 
representing sheaf.  

We end the paper with the following two remarks which are often
considered in the study of this course.

{\bf Remark~1.} We would like to point out a special case
where $\mathcal M$ is a twisted structure sheaf $\mathcal O(-m)$. This
has drawn the attention of those who attempted to solve Exercise~19.18
in \cite{E}. For any $m \in \mathbb Z$, we obtained
\begin{equation}\label{eqpoly}
P_{\mathcal O(-m)}(t) = \sum_{\ell=0}^m (-1)^{\ell} \binom{m}{\ell}
P_{S_{\ell}}(t) 
\end{equation} 
by an inductive computation on binomial coefficient functions
\begin{equation}\label{eqbinom}
\binom{t+d-m}{d}=\binom{t+d-(m-1)}{d} - \binom{t+(d-1)-(m-1)}{d-1}.
\end{equation}
We use the convention that $\binom{a}{b}=0$ if $a<b$.
Part~B in Proposition~\ref{chernhilbert} can be reduced to a problem
asking the following congruence
\begin{equation}\label{twisted} 
C_{\mathcal O(-m)}(x)=1-mx \equiv 
\prod_{\ell=0}^m (C_{S_{\ell}})^{(-1)^{\ell}\binom{m}{\ell}} ,
\pmod{x^{d+1}}. 
\end{equation} 

Since $C_{S_{\ell}}$ can be computed by the Koszul complex as it is
shown in
(\ref{eqchern}), a naive attempt on proving the congruence in the
above (\ref{twisted}) is 
to compute the coefficients of each $x^i$ on the right hand side
and to show that the coefficients of higher terms are zero. However,
the coefficient of a general term $x^i$ in terms of the binomial
coefficients is rather complicated. It is not clear how these terms
vanish for $i \geq 2$ if they are treated as combinatorial formulae.

The result in the previous section provides the following intuition from a 
different perspective: (\ref{twisted}) follows from the fact
that $\sigma$ is a group isomorphism; precisely, 
\[ \begin{array}{lcl}
 (C_{\mathcal O(-m)}(x),1)=(1-mx,1)& = & 
 \sigma^{-1}(P_{\mathcal O(-m)}(t)) = 
 \sigma^{-1}(\sum_{\ell=0}^m (-1)^{\ell}\binom{m}{\ell}P_{S_{\ell}}(t)) \\
 & = & (\prod_{\ell
    =0}^m (C_{S_{\ell}}(x))^{(-1)^{\ell}\binom{m}{\ell}}, 1). 
   \end{array} \]
Hence, $C_{\mathcal O(-m)}(x)$ and $\prod_{\ell
    =0}^m (C_{S_{\ell}}(x))^{(-1)^{\ell}\binom{m}{\ell}}$ are equal in
  the groups $\mathcal A_0$. 
Although (\ref{eqbinom}) is a combinatorial property, it is also the
 relations of Hilbert polynomials of the sheaves in the following short
 exact sequence
\begin{equation}\label{hyperplane}
 0 \longrightarrow \mathcal O(-m) \stackrel{\mathcal H}{\longrightarrow}
 \mathcal O(-m+1) \longrightarrow \mathcal O_{\mathcal
  H}(-m+1) \longrightarrow 0 , \end{equation}
where $\mathcal H$ is a hyperplane in $\mathbb P^d$. 
Since the Chern polynomials depend on the ambient scheme, a similar
inductive decomposition as it is for $P_{\mathcal M}(t)$ in (\ref{eqpoly}) 
does not hold for Chern polynomials. However, (\ref{hyperplane}) induces 
an identity on the cycles $\alpha = [\mathcal O(-m)] =\sum_{\ell=0}^m
(-1)^{\ell}\binom{m}{\ell} [S_{\ell}] $ in $\K(\mathbb P^d)$ which implies
the corresponding (\ref{twisted}) in $\mathcal A$.

{\bf Remark~2.} A more fundamental correspondence of the
Chern and Hilbert polynomials should be pointed out. Let $a_i$ denote
the coefficient of $x^i$ in the Taylor expansion of
$(\frac{x}{1-e^{-x}})^{d+1}$. (\ref{eqPC}) can be written explicitly as
\begin{equation}\label{Pcoeff} \begin{array}[t]{ll} 
P_M(t)= & \frac{1}{d!}a_0r t^d + \frac{1}{(d-1)!}(a_0 \ch_1+a_1r)t^{d-1}
 + \\ \\ & \cdots 
(a_0\ch_d+a_1\ch_{d-1}+\cdots + a_{d-1}\ch_1 + a_d r) , \end{array}
\end{equation}
in which we abbreviate $\ch_i(\mathcal M)$ by $\ch_i$. 
Replacing the above $\ch_i$ by the proper terms in (\ref{chernchar}),
the coefficients of $P_{\mathcal M}$ can be expressed in terms of Chern
classes. Conversely, if $P_{\mathcal M}$ is known; that is, the
coefficients of $P_{\mathcal M}(t)$ are determined, then the Chern
classes can be solved inductively using the above (\ref{Pcoeff}).
However, the computation for $a_i$ is tedious.
Part~B in Proposition~\ref{chernhilbert} avoids such lengthy
computation.

\bigskip

\begin{center} {\sc Acknowledgment.} \end{center} 
\noindent {\em The author would like to express her gratitude to David
Eisenbud and Paul Roberts for introducing to her this interesting
problem and for many valuable discussions. To Donu Arapura, Vesselin
Gasharov, Kazuhiko Kurano and Kenji Matsuki for their important
suggestions in the proof. She also thanks the referee on the revision
of the paper in the current form. Part of the paper was written while
the author was visiting I-Chiau Huang at Academia Sinica in
Taiwan. Their hospitality is greatly appreciated.}

\bigskip
{\small 
{\em Department of Mathematics, Purdue University, West Lafayette, 
 IN~47907-1395} \\
\indent E-mail address: {\em chan@math.purdue.edu} \\
\indent Current address: \\
\indent {\em Department of Mathematics, University of Arkansas, Fayetteville, 
  AR~72701.}}

\end{document}